\documentclass{amsart}

\usepackage{amsfonts}
\usepackage{amsmath}
\usepackage{amsthm}
\usepackage{amssymb}
\usepackage[english]{babel}
\usepackage{graphics}
\usepackage{latexsym}
\usepackage{longtable} \usepackage{mathrsfs}
\usepackage{graphicx}
\usepackage{wrapfig} 
\usepackage[pdftex,colorlinks=true]{hyperref}

\newtheorem{theorem}{Theorem}
\newtheorem{corollary}[theorem]{Corollary}
\newtheorem{lemma}[theorem]{Lemma}
\newtheorem{proposition}[theorem]{Proposition}

\newtheorem{claim}[theorem]{Claim}
\newtheorem{example}[theorem]{Example}

\theoremstyle{definition}
\newtheorem{definition}[theorem]{Definition}
\newtheorem{remark}[theorem]{Remark}

\renewcommand{\S}{\mathcal{S}}

\renewcommand{\H}{\textrm{H}}
\newcommand{\R}{\mathbb{R}}
\newcommand{\N}{\mathbb{N}}
\newcommand{\mS}{\mathbb{S}}
\newcommand{\mB}{\mathbb{B}}

\newcommand{\noi}{\noindent}
\newcommand{\ms}{\medskip}

\newcommand{\al}{\alpha}
\newcommand{\be}{\beta}
\newcommand{\ga}{\gamma}
\newcommand{\de}{\delta}
\newcommand{\De}{\Delta}
\newcommand{\e}{\varepsilon}

\newcommand{\Om}{\Omega}

\newcommand{\larrow}{\longrightarrow}
\newcommand{\ot}{\otimes}

\newcommand{\ri}{\rightarrow}

\newcommand{\p}{\partial}
\newcommand{\sub}{\subseteq}
\newcommand{\set}{\setminus}

\newcommand{\rk}{\textrm{rk}}

\newcommand{\Lip}{\textrm{Lip}}
\newcommand{\tr}{\textrm{tr}}
\newcommand{\sgn}{\textrm{sgn}}
\newcommand{\ess}{\textrm{ess}}

\newcommand{\dist}{\textrm{dist}}

\newcommand{\Div}{\textrm{Div}}

\newcommand{\bt}{\begin{theorem}}\newcommand{\et}{\end{theorem}}
\newcommand{\bd}{\begin{definition}}\newcommand{\ed}{\end{definition}}
\newcommand{\bl}{\begin{lemma}}\newcommand{\el}{\end{lemma}}
\newcommand{\beq}{\begin{equation}}\newcommand{\eeq}{\end{equation}}
\newcommand{\bc}{\begin{claim}}\newcommand{\ec}{\end{claim}}
\newcommand{\bex}{\begin{example}}\newcommand{\eex}{\end{example}}
\newcommand{\bcor}{\begin{corollary}}\newcommand{\ecor}{\end{corollary}}
\newcommand{\bp}{\begin{proof}}\newcommand{\ep}{\end{proof}}

\newcommand{\BPL}{\medskip \noindent \textbf{Proof of Lemma} }

\newcommand{\BPP}{\medskip \noindent \textbf{Proof of Proposition} }
\newcommand{\BPT}{\medskip \noindent \textbf{Proof of Theorem} }

\numberwithin{equation}{section}
\numberwithin{theorem}{section}

\begin{document}

\title{$\infty$-Minimal Submanifolds}

\author{\textsl{Nikolaos I. Katzourakis}}
\address{BCAM - Basque Center for Applied Mathematics, Alameda de Mazarredo 14, E-48009, Bilbao, Spain}
\email{nkatzourakis@bcamath.org}

\subjclass[2010]{Primary 35J47, 35J62, 53C24; Secondary 49J99}

\date{}


\keywords{$\infty$-Harmonic maps, Vector-valued Calculus of Variations in $L^\infty$,  Vector-valued Optimal Lipschitz Extensions, Quasi-Conformal maps, Aronsson PDE, Rigidity.}

\begin{abstract} We identify the Variational Principle governing $\infty$-Harmonic maps $u : \Om \sub \R^n \larrow \R^N$, that is solutions to the $\infty$-Laplacian
\[  \label{1}
\De_\infty u \ :=\ \Big(Du \ot Du \, +\, |Du|^2 [Du]^\bot \! \ot I \Big) : D^2 u\ = \ 0. \tag{1}
\]
System \eqref{1} was first derived in the limit of the $p$-Laplacian as $p\ri \infty$ in \cite{K2} and was recently studied in \cite{K3}. Here we show that \eqref{1} is the ``Euler-Lagrange PDE'' of vector-valued Calculus of Variations in $L^\infty$ for the functional 
\[ \label{2}
\|Du\|_{L^\infty(\Om)}\ = \ \underset{\Om}{\ess\,\sup} \,|Du|.     \tag{2}
\]
We introduce the notion of \emph{$\infty$-Minimal Maps}, whch are Rank-One Absolute Minimals of \eqref{2} with \emph{``$\infty$-Minimal Area''} of the submanifold $u(\Om) \sub \R^N$ and prove they solve \eqref{1}. The converse is true for immersions. We also establish a maximum principle for $|Du|$ for solutions to \eqref{1}. We further characterize minimal surfaces of $\R^3$ as those locally parameterizable by isothermal immersions with $\infty$-Minimal area and show that isothermal $\infty$-Harmonic maps are rigid.
\end{abstract}

\maketitle

\section{Introduction} \label{section1}

 In this paper we are interested in the variational structure of \emph{$\infty$-Harmonic maps},  that is of solutions $u : \Om \sub \R^n \larrow \R^N$, $n,N\geq 2$, to the PDE system
\beq  \label{1.1}
\De_\infty u \ :=\ \Big(Du \ot Du \, +\, |Du|^2 [Du]^\bot \! \ot I \Big) : D^2 u\ = \ 0. 
\eeq
Here $[Du(x)]^\bot$ is the projection on the nullspace of the transpose of the gradient matrix $Du(x)^\top : \R^N \larrow \R^n$ and $|Du|^2=\tr(Du^\top Du)$ is the Euclidean norm on $\R^N \ot \R^n$ (for details see Preliminaries \ref{Preliminaries}). In index form,  \eqref{1.1} reads
\beq  
 D_i u_\al  D_j u_\be D_{ij}^2 u_\be \, +\, |Du|^2 [Du]_{\al \be}^\bot D^2_{ii} u_\be\ = \ 0
\eeq
with triple summation in $1\leq i,j\leq n$ and $1\leq \be \leq N$. System \eqref{1.1} is a quasilinear degenerate elliptic system in non-divergence form which arises in the limit of the $p$-Laplace system $\De_p u = \Div \big(|Du|^{p-2}Du\big)=0$ as $p\ri \infty$. It was first derived by the author in \cite{K2} and was studied in the very recent work \cite{K3}. The special case of the scalar $\infty$-Laplace PDE for $N=1$ reads
\beq \label{1.2}
\De_\infty  u \ =\ Du \ot Du :D^2u\ = \ 0 
\eeq 
and has a long history. In this case the coefficient $|Du|^2[Du]^\bot$ of \eqref{1.1} vanishes identically and the same holds for submersions in general. Equation \eqref{1.2} was derived in the limit of the $p$-Laplacian as $p\ri \infty$ in the '60s by Aronsson and was first studied in \cite{A3, A4}. It has been extensively studied ever since, in the last 20 years in the context of Viscosity Solutions (see for example Crandall \cite{C}, Barron, Evans, Jensen \cite{BEJ} and references therein). A major difficulty in its study is its degeneracy and the emergence of  singular solutions (see e.g.\ \cite{K1}).

Aronsson derived \eqref{1.2} in the limit of the Euler-Lagrange equation of the $p$-Dirichlet functional, or equivalently of the $L^p$-norm of the gradient $\|Du\|_{L^p(\Om)}$. He observed that at least in a formal level $\De_p u\ri \De_\infty u$ and  $\|Du\|_{L^p(\Om)}\ri \|Du\|_{L^\infty(\Om)}$ both as $p\ri \infty$, but it was not a priori clear that the following rectagle ``commutes" 
\begin{align} \label{1.3}
&\|Du\|_{L^p(\Om)}   \ \ \ \  \larrow\ \ \ \ \De_p u =0  \nonumber\\
& \ \ \ \downarrow\ p \ri \infty\ \ \ \ \ \ \ \ \ \ \  \ \ \ \ \ \ \downarrow\ p \ri \infty  \\
&\|Du\|_{L^\infty(\Om)}   \ \ \ \ \dashrightarrow \ \ \ \ \De_\infty u=0   \nonumber
\end{align}
so that \eqref{1.2} has a variational structure with respect to the \emph{supremal functional}
\beq \label{1.4}
\|Du\|_{L^\infty(\Om)}\ = \ \underset{\Om}{\ess\,\sup} \,|Du|,
\eeq
in the sense that \eqref{1.2} is the ``Euler-Lagrange PDE'' of Calculus of Variations in $L^\infty$ for the model functional \eqref{1.4}. This turned out to be the case and inspired by his earlier work \cite{A1,A2} he identified the appropriate variational notion, that of \emph{Absolute Minimals} for \eqref{1.4}, which alllows to connect \eqref{1.4} with \eqref{1.2}. The subtle point is that \eqref{1.4} is \emph{nonlocal}, in the sense that with respect to the $\Om$ argument \eqref{1.4} is not a measure. This implies that minimizers over a domain with fixed boundary values are not local minimizers over subdomains and the direct method of Calculus of Variations when applied to \eqref{1.4} does not produce PDE solutions of \eqref{1.2}. Absolute Minimals is nothing but local minimizers of \eqref{1.4}, but locality is built into the minimality notion:
\beq \label{1.5}
\left.
\begin{array}{l}
D \subset \subset \Om, \\
g\in W^{1,\infty}_0(D)
\end{array}
\right\} \ \ \Longrightarrow \ \
\big\|Du\big\|_{L^\infty(D)}\ \leq \ \big\|D(u+g)\big\|_{L^\infty(D)}.
\eeq
Aronsson established the equivalence between Absolute Minimals satisfying \eqref{1.5} and solutions to \eqref{1.2}, namely $\infty$-Harmonic functions, in the smooth setting. This result was later extended to general viscosity solutions of \eqref{1.2} (see \cite{C}).

In the full vector case of \eqref{1.1}, even more intriguing phenomena occur, studied in the case of smooth solutions in \cite{K2,K3}. Except for the emergence of ``singular solutions'' to \eqref{1.1}, a further difficulty not present in the scalar case is that \emph{\eqref{1.1} has discontinuous coefficients} even for $C^\infty$ solutions. There exist smooth $\infty$-Harmonic maps whose rank of the gradient is not constant: such an example on $\R^2$ is given by $u(x,y) = e^{ix}-e^{iy}$. This $u$ is $\infty$-Harmonic near the origin and has $\rk(Du)=1$ on the diagonal, but it has $\rk(Du)=2$ otherwise and hence the projection $[Du]^\bot$ is discontinuous. In general, \emph{$\infty$-Harmonic maps present a phase separation},  studied for $n=2\leq N$ in \cite{K3}. On each phase the dimension of the tangent space is constant and these phases are separated by \emph{interfaces} whereon the rank of $Du$ ``jumps'' and $[Du]^\bot$ gets discontinuous. 
\[
\underset{\text{Figure 1.}}{\includegraphics[scale=0.2]{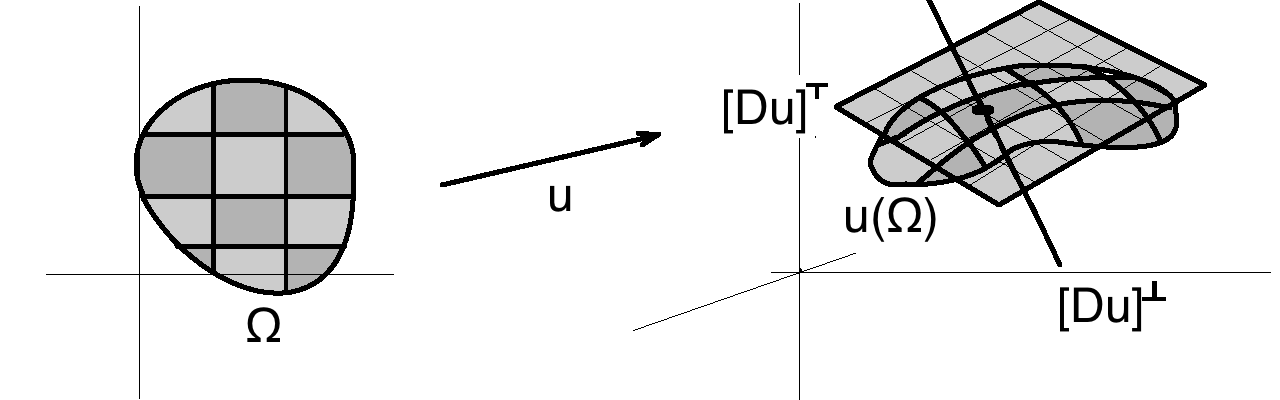}}
\]
On a phase, we interpret \eqref{1.1} as decoupling to the \emph{tangential} system $Du D\big(\frac{1}{2}|Du|^2\big)$ $=0$ in the tangent bundle $[Du]^\top$ and the \emph{normal} system $|Du|^2 [Du]^\bot \De u=0$ in the normal bundle $[Du]^\bot$.   Interestingly, discontinuous coefficients is a genuine vectorial phenomenon of general maps and does not arise when either $n=1$ or $N=1$. In particular, when $n=1$  all $\infty$-Harmonic curves are \emph{affine} and for $u :\Om \sub \R \larrow \R^N$, \eqref{1.1} reduces to 
\begin{align} \label{1.6}
\De_\infty u \ & =\ (u' \ot u') u''\, +\, |u'|^2\Big(I - \frac{u'}{|u'|}\ot \frac{u'}{|u'|}\Big)u''\ =\ |u'|^2u''.
\end{align}

In this paper we identify the appropriate variational notion for the model functional \eqref{1.4} of vector-valued Calculus of Variations in $L^\infty$ which characterizes system \eqref{1.1} and also consider some related questions. In the case $N>1$, we equip $\R^N \ot \R^n$ with the Euclidean norm. In \cite{K2} we established that Aronsson's notion of Absolute Minimals adapted to the vector case indeed leads to solutions of the tangential system $Du\ot Du:D^2u=0$, but the question of how to describe variationally the full system \eqref{1.1} remained open. We also showed that $Du\ot Du:D^2u=0$ is not sufficient for Absolute Minimality. 

Herein we settle these problems. In Definition \ref{def1} we introduce the variational notion of  \emph{$\infty$-Minimal Maps $u : \Om \sub \R^n \larrow \R^N$}. An $\infty$-Minimal map is a weak version of Absolute Minimal of \eqref{1.4} with respect to essentially scalar local variations with zero boundary values which we call \emph{Rank-One Absolute Minimal} (Definition \ref{def1} $(i)$) coupled by a notion of \emph{``$\infty$-Minimal Area'' of the submanifold $u(\Om) \sub \R^N$}  (Definition \ref{def1} $(ii)$). The latter means minimality  for \eqref{1.4} with respect to local variations normal to $u(\Om)$ with free boundary values. In order for these conditions to be made rigorous and precise, we restrict ourselves to the case of smooth maps \emph{of full rank}, that is when $\rk(Du)=\min\{n,N\}$. This class consists of immersions, submersions and local diffeomorphisms. With a little extra effort we could consider smooth maps $u$ where the rank of $Du$ is any piecewise constant function on sets with nonempty interior, but the difficulty of discontinuous coefficients of \eqref{1.1} comes into play and we can not go much further without an appropriate ``weak'' theory of nondifferentiable solutions of system \eqref{1.1}.

In Theorem \ref{th1} we prove that $\infty$-Minimal Maps are solutions to \eqref{1.1}. The converse is true for immersions and \eqref{1.1} is both necessary and sufficient for the variational problem in this class. Rather surprisingly, for immersions $\infty$-Minimality of the area is also equivalent to a relevant notion of \emph{$p$-Minimal Area of $u(\Om)$} for all $p\in [2,\infty)$, where normal variations are considered for the $L^p$ norm of the gradient. Moreover, in Proposition \ref{pr1} we establish a maximum and a minimum principle for $|Du|$ of solutions $u$ to \eqref{1} with full rank, by employing an improved version of the gradient flow introduced in \cite{K2}, which bears the property of the scalar case that (projections of) \emph{images of trajectories $t\mapsto \xi^\top u(\ga(t))$ are affine}. 

The conditions of $p$- and $\infty$-Minimal area of $u(\Om)$ are definitely reminishent to that of Minimal Surfaces. In the case of the latter, what we consider is normal variations of the Area of the surface, which is the integral of the \emph{Jacobian}. Interestingly, in the class of \emph{conformal} maps $u : \Om \sub \R^2 \larrow \R^3$, the quantity $[Du]^\bot\De u$ is proportional to the \emph{mean curvature vector $\H$ of $u(\Om)$}, while the Area coincides with the Dirichlet functional (Lemma \ref{l4}, Corollary \ref{cor2}). These observations allow us to characterize minimal surfaces $\S$ of $\R^3$ as those locally parameterizable by isothermal normally $\infty$-Harmonic maps and surfaces which are locally isometric to $\R^2$ as those locally parameterizable by isothermal tangentially $\infty$-Harmonic maps (Theorem \ref{th2}). As a corollary, we deduce a \emph{rigidity result}: isothermal $\infty$-Harmonic maps have affine image (Corollary \ref{cor3}).

We conclude this introduction by recalling some very recent important vectorial results related to \eqref{1.1} and \eqref{1.4}.  Ou, Troutman and Wilhelm in \cite{OTW} and Wang and Ou in \cite{WO} studied Riemannian variants of tangentially $\infty$-Harmonic maps which solve only the tangential part of \eqref{1.1}. Sheffield and Smart in \cite{SS} used the nonsmooth operator norm on $\R^N \ot \R^n$ and derived a singular variant of \eqref{1.1} connected to $\ess\, \sup_\Om \|Du\|$ for a norm different than the Euclidean, which governs \emph{optimal Lipschitz extensions of maps}. The authors use this norm because they need the coincidence of $\|Du\|_{L^\infty(\Om)}$ with the Lipschitz constant $\Lip(u,\Om)$, which fails for the Euclidean norm $|Du|$ on $\R^N \ot \R^n$. They introduced the optimality notion of \emph{tightness} for Lipschitz extensions and characterized smooth solutions of their version of $\De_\infty$ as tight maps. Capogna and Raich in \cite{CR} used the supremal functional $\ess\,\sup_\Om  {|Du|^n}/{\textrm{det}(Du)}$ defined for local diffeomorphisms $u : \Om \sub \R^n \larrow \R^n$ and developed an $L^\infty$ variational approach to extremal \emph{Quasi-Conformal maps}. They derived a variant of \eqref{1.1}, for which the normal term vanishes identically and studied smooth extremal Quasi-Conformal maps as solutions of an Aronsson system.  Their results have very recently been advanced by the author in \cite{K4}.

\subsection{Preliminaries.} \label{Preliminaries} Throughout this paper we reserve $n,N \in \N$ for the dimensions of Euclidean spaces and $\mS^{N-1}$ denotes the unit sphere of $\R^N$. Greek indices $\al, \be, \ga,... $ run from $1$ to $N$ and Latin $i,j,k,...$ form $1$ to $n$. The summation convention will always be employed in repeated indices in a product. Vectors are always viewed as columns. Hence, for $a,b\in \R^n$, $a^\top b$ is their inner product and $ab^\top$ equals $a \ot b$.  If $V$ is a vector space, then $\mS(V)$ denotes the symmetric linear maps $T : V \larrow V$ for which $T=T^\top$. If $u=u_\al e_\al :  \Om \sub \R^n \larrow \R^N$ is in $C^2(\Om)^N$, the gradient matrix $Du$ is viewed as $D_i u_\al e_\al \ot e_i : \Om \larrow \R^N \ot \R^n$ and the Hessian tensor $D^2u$ as $D^2_{ij} u_\al e_\al \ot e_i \ot e_j: \Om \larrow \R^N \ot \mS(\R^n)$. The Euclidean (Frobenious) norm on $\R^N \ot \R^n$ is $|P|=(P_{\al i}P_{\al i})^{\frac{1}{2}} = (\tr (P^\top P))^{\frac{1}{2}}$.  We also introduce the following \emph{contraction operation} for tensors which extends the Euclidean inner product $P:Q=\tr(P^\top Q)=P_{\al i}Q_{\al i}$ of $\R^N \ot \R^n$. Let ``$\ot^{(r)}$'' denote the $r$-fold tensor product. If $S\in \ot^{(q)}\R^N \ot^{(s)} \R^n$,  $T \in \ot^{(p)}\R^N \ot^{(s)} \R^n$ and $q\geq p$, we define a tensor $S:T$ in $\ot^{(q-p)} \R^N$ by 
\beq \label{1.24}
S:T \ :=\ \big(S_{\al_q ...\al_p... \al_1 \, i_s ... i_1}  T_{\al_{p}  ... \al_1 \, i_s ... i_1}  \big) \, e_{\al_q} \ot ... \ot e_{\al_{p+1}}. 
\eeq
For example, for $s=q=2$ and $p=1$, the tensor $S:T$ of \eqref{1.24} is a vector with components $S_{\al \be i j}T_{\be ij}$ with free index $\al$ and  the indices $\be,i,j$ are contracted. In particular, in view of \eqref{1.24}, the 2nd order linear system 
\beq
A_{\al i \be j}D^2_{ij}u_\be \, +\, B_{\al \ga k} D_ku_\ga + C_{\al \de} u_\de\, =\, f_\al ,
\eeq
can be compactly written as $A$:$D^2u + B$:$Du+Cu=f$, where the meaning of ``$:$" in the respective dimensions is made clear by the context. Let now $P : \R^n \larrow \R^N$ be linear map. Upon identifying linear subspaces with orthogonal projections on them (with respect to the standard inner product), we split $\R^N=[P]^\top \oplus [P]^\bot$ where $[P]^\top$ and $[P]^\bot$ denote range of $P$ and nullspace of $P^\top$ respectively. Hence, if $\xi \in \mS^{N-1}$, then $[\xi]^\bot$ is (the projection on) the normal hyperplane $I-\xi \ot \xi$. Let now $u : \Om \sub \R^n \larrow \R^N$ be a map in $C^1(\Om)^N$. Generally, the rank of $Du$ satisfies $\rk(Du) \leq \min\{n,N\}$. We will call $u$ a \emph{Full-Rank Map} if $\rk(Du)=\min\{n,N\}$ on $\Om$, that is when $\rk(Du)$ achieves the maximum possible value everywhere on $\Om$. If $n\leq N$ then $u$ is an \emph{immersion} and if $n\geq N$ then $u$ is a \emph{submersion}. If both happen and $n=N$, then $u$ is a \emph{local diffeomorphism}. For immersions, the \emph{Jacobian} $Ju$ is the square root of the determinant of the induced from $\R^N$ Riemannian metric on $u(\Om)$, that is $Ju:= \sqrt{\det(Du^\top Du)}$. The map $u$ is \emph{Conformal} when there is $f\in C^1(\Om)$ such that $Du^\top Du =f^2 I$ on $\Om$, that is $D_i u_\al D_j u_\al=f^2\de_{ij}$. If $n=2$, $N=3$ and $f\neq0$, conformal immersions are called isothermal parametrizations of the surface $u(\Om)\sub \R^3$. Given a full-rank map $u$, we will identify the pull back of the tangent bundle of $u(\Om)$ to $\Om$ with the projection $[Du]^\top$ and its orthogonal complement with the projection $[Du]^\bot$. We will denote the set of tangent vector fields along $u$ by $\Gamma([Du]^\top)$ and  the set of normal vector fields along $u$ by $\Gamma([Du]^\bot)$. Obviously, if $u$ is a submersion, then $\Gamma([Du]^\bot)$ contains only the zero vector field.

\section{Variational Structure of $\infty$-Harmonic Maps.} \label{section2}

We begin by introducing a minimality notion of vector-valued Calculus of Variations in $L^\infty$ for the supremal functional
\beq \label{2.1}
\|Du\|_{L^\infty(\Om)} \ = \ \underset{\Om}{\ess\, \sup}\, |Du|,
\eeq 
where $|Du|$ is the Euclidean norm on $\R^N \ot \R^n$. 

\begin{definition}\label{def1} Let $u : \Om \sub \R^n \larrow \R^N$ be a map in $C^1(\Om)^N$.

\ms \noi (i) The map $u$ is called \emph{Rank-One Absolute Minimal} on $\Om$ when for all compactly contained subdomains $D$ of $\Om$, all functions $g$ on $D$ vanishing on $\p D$ and all directions $\xi$, $u$ is a minimizer on $D$ with respect to essentially scalar variations $u+g\xi$:
\beq \label{2.2}
\left.
\begin{array}{l}
D \subset \subset \Om, \\
g\in C^1_0(D), \\
\xi \in \mS^{N-1}
\end{array}
\right\} \ \ \Longrightarrow \ \
\big\|Du\big\|_{L^\infty(D)}\ \leq \ \big\|D(u+g\xi)\big\|_{L^\infty(D)}.
\eeq
\[
\underset{\text{Figure 2.}}{\includegraphics[scale=0.18]{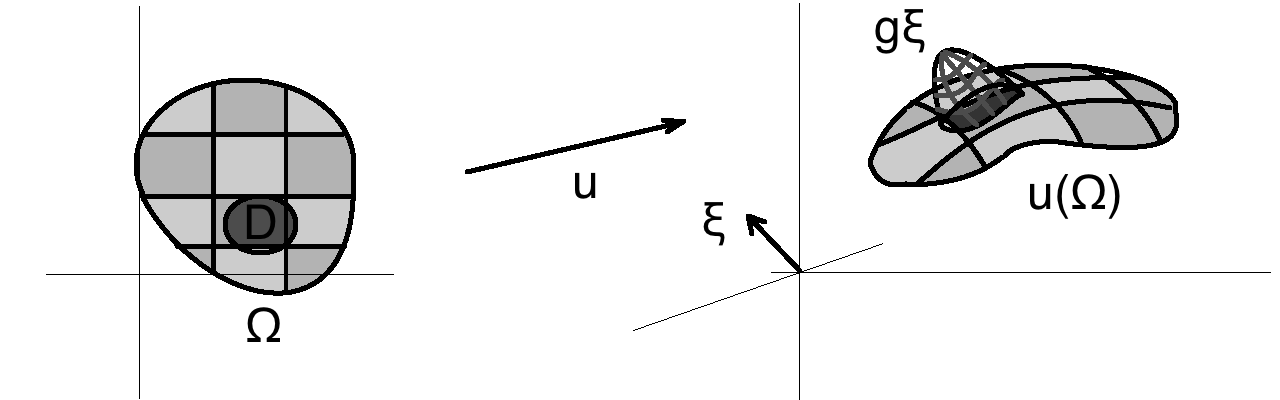}}
\]
 \noi (ii) Suppose $u$ is an immersion. We say that \emph{$u(\Om)$ has $\infty$-Minimal Area} when for all compactly contained subdomains $D$, all functions $h$ on $\bar{D}$ (not only vanishing on $\p D$) and all normal vector fields $\nu$, $u$ is a minimizer on $D$ with respect to normal free variations $u+h\nu$:
\beq \label{2.3}
\left.
\begin{array}{l}
D \subset \subset \Om, \\
h\in C^1(\bar{D}), \\
\nu \in \Gamma([Du]^\bot)
\end{array}
\right\} \ \ \Longrightarrow \ \
\big\|Du\big\|_{L^\infty(D)}\ \leq \ \big\|D(u+h\nu)\big\|_{L^\infty(D)}.
\eeq
\[
\underset{\text{Figure 3.}}{\includegraphics[scale=0.18]{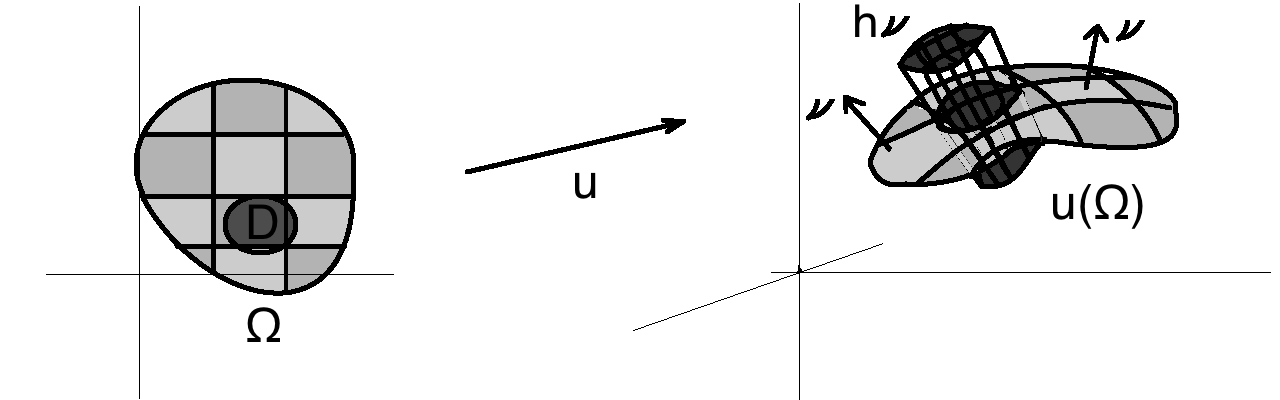}}
\]
Similarly, if \eqref{2.3} holds with the $L^p$ norm in the place of the $L^\infty$ norm, we will say that the image $u(\Om)$ of the immersion has ``\emph{$p$-Minimal Area}''.

\ms

 \noi (iii) Suppose $u$ is a Full-Rank map, that is $\rk(Du)=\min\{n,N\}$ on $\Om$. Then, we call $u$ an \emph{$\infty$-Minimal Map with respect to functional \eqref{2.1}} when $u$ is a Rank-One Absolute Minimal on $\Om$ and $u(\Om)$ has $\infty$-Minimal Area.
\end{definition}

Evidently, condition $(ii)$ of Definition \ref{def1} is empty for submersions and in particular in the scalar case $N=1$, since in $0$-codimension we have $\rk(Du)=N\leq n$ and hence $[Du]^\bot=\{0\}$ in this case.

\bt[Variational Structure of $\infty$-Laplacian] \label{th1}  Let $u : \Om \sub \R^n \larrow \R^N$ be a  map in $C^2(\Om)^N$. Then:

\ms
\noi (i) If $u$ is an $\infty$-Minimal Map with respect to functional $\|Du\|_{L^\infty(\Om)}$, it follows that $u$ is $\infty$-Harmonic on $\Om$ and solves the system
\beq  \label{2.4}
\De_\infty u \ =\ \Big(Du \ot Du \, +\, |Du|^2 [Du]^\bot \! \ot I \Big) : D^2 u\ = \ 0. 
\eeq
If $u$ is an immersion, the converse is true as well and $\infty$-Harmonicity implies $\infty$-Minimality.
In particular, the following assertions hold for the tangential and the nornal part separately:

\ms \noi (ii) If $u$ is a Rank-One Absolute Minimal on $\Om$, then $u$ is tangentially $\infty$-Harmonic  on $\Om$ and solves $Du \ot Du :D^2u= 0$. The converse is true if $u$ is an immersion.
\ms

\noi (iii) Suppose $u$ is an immersion. Then, $u(\Om)$ has $\infty$-Minimal Area if and only if $u$ is normally $\infty$-Harmonic  on $\Om$ and solves $|Du|^2 [Du]^\bot  \De u= 0$. 
\et

\noi We note that for immersions the system $|Du|^2 [Du]^\bot  \De u= 0$ is equivalent to $[Du]^\bot  \De u= 0$, but we keep the positive function $|Du|^2$ because for ``singular solutions'' these systems generally are not equivalent. 

The proof of Theorem \ref{th1} is split in four lemmas. The first one below is implied by Theorem 2.1 of \cite{K2}, but for the sake of completeness we provide a sharper simplified proof. 

\bl \label{l1} Let $u : \Om \sub \R^n \larrow \R^N$ be in $C^2(\Om)^N$. If $u$ is a Rank-One Absolute Minimal, then $u$ is tangentially $\infty$-Harmonic and solves $Du\ot Du :D^2u=0$ on $\Om$.
\el

\BPL \ref{l1}. Fix $x\in \Om$, $0<\e < \dist(x,\p \Om)$, $0<\de<1$ and $\xi \in \mS^{N-1}$. Choose $D:=\mB_\e(x)$, $g(z):=\frac{\de}{2}\big(\e^2-|z-x|^2\big) \in C^1_0(D)$ and set $w:=u+g\xi$. Then, by Taylor expansions of $|Du|^2$ and $|Dw|^2$ at $x$ we have
\begin{align} \label{2.7}
|Du(z)|^2\ = \ |Du(x)|^2\ +\ D\big(|Du|^2\big)(x)^\top(z-x)\ + \ o(|z-x|),
\end{align}
as $z\ri x$, and also by using that $D^2g=-\de I$ and $Dg(x)=0$ we have
\begin{align} \label{2.8}
|Dw(z)|^2\ & = \ |Du(x)+ \xi \ot Dg(x)|^2\ +\ D\big(|Du+ \xi \ot Dg|^2\big)(x)^\top(z-x) \nonumber\\
& \ \ \ \  + \ o(|z-x|)\nonumber\\
             & = \ |Du(x)|^2\ +\ 2Du(x)^\top \big(D^2u(x)-\de \xi\ot I \big)(z-x) \\\
     &\ \ \ \  + \ o(|z-x|)\nonumber\\
 & = \ |Du(x)|^2\ + \ \Big(D\big(|Du|^2\big)(x)^\top  -\, 2\de \xi^\top Du(x)\Big) (z-x) \nonumber\\
&\ \ \ \  + \ o(|z-x|),\nonumber
\end{align}
as $z\ri x$. By \eqref{2.7} we have the estimate
\begin{align} \label{2.9}
\|Du\|_{L^\infty(\mB_\e(x))}^2\ & \geq  \ |Du(x)|^2 \ +\ \max_{\{|z-x|\leq \e\}}\Big\{D\big(|Du|^2\big)(x)^\top(z-x)\Big\}\  \nonumber\\
&\ \ \ \ \ + \ o(\e)\\
& = \ |Du(x)|^2\ +\ \e\big|D\big(|Du|^2\big)(x)\big|\ + \ o(\e),\nonumber
\end{align}
as $\e \ri 0$, and also by \eqref{2.8} we have
\begin{align} \label{2.10}
\|Dw\|_{L^\infty(\mB_\e(x))}^2\ & \leq  \ |Du(x)|^2 \ + \max_{\{|z-x|\leq \e\}}\Big\{\big(D\big(|Du|^2\big)(x)^\top \! - 2\de \xi^\top Du(x)\big)(z-x)\Big\} \nonumber\\
&\ \ \ \  + \ o(\e) \\
& = \ |Du(x)|^2\ +\ \e\big|D\big(|Du|^2\big)(x) - \, 2\de \xi^\top Du(x)\big|\ + \ o(\e) , \nonumber
\end{align}
as $\e \ri 0$. Then, since $u$ is Rank-One Absolute Minimal on $\Om$,  inequalities \eqref{2.9} and  \eqref{2.10}  imply
\begin{align} \label{2.11}
0\ & \leq \ \|Dw\|_{L^\infty(\mB_\e(x))}^2 \, - \,  \|Du\|_{L^\infty(\mB_\e(x))}^2 \nonumber\\
&\leq \  \e\Big(\big|D\big(|Du|^2\big)(x) -\, 2\de \xi^\top Du(x)\big|\ -\ \big|D\big(|Du|^2\big)(x)\big|\Big)\ + \ o(\e),
\end{align}
as $\e \ri 0$. If $D\big(|Du|^2\big)(x)=0$, we obtain $ (Du \ot Du : D^2u)(x) = 0$ as desired. If $D\big(|Du|^2\big)(x)\neq0$, then Taylor expansion of $p \mapsto  \big|D\big(|Du|^2\big)(x) +\, p\big| - \big|D\big(|Du|^2\big)(x)\big|$ at $p_0=0$ and evaluated at $p=-\, 2\de \xi^\top Du(x)$, \eqref{2.11} implies after letting $\e \ri 0$ that
\beq \label{2.12a}
0\ \leq  \ - 2\de \, \xi^\top Du(x)\left(\frac{D\big(|Du|^2\big)(x)}{\big|D\big(|Du|^2\big)(x)\big|} \right) \ + \ o(\de).
\eeq
By letting $\de \ri 0$ in \eqref{2.12a} we obtain $ \big( \xi^\top Du \ot Du : D^2u\big)(x) \geq 0$ and since $\xi$ is arbitrary we get $ (Du \ot Du : D^2u)(x) = 0$ for any $x\in \Om$. The lemma follows.               \qed

\ms
Now we consider the converse of Lemma \ref{l1}, that is the sufficiency of the tangential part of the $\infty$-Laplacian for Rank-One Absolute Minimality. Example 3.3 in \cite{K2} shows that $Du \ot Du :D^2u=0$ does \emph{not imply} the stronger condition of Absolute Minimality with respect to arbitrary vectorial variations. Lemma \ref{l1a} below is valid only for the weaker rank-one condition of essentially scalar variations.

\begin{lemma} \label{l1a} Let $u : \Om \sub \R^n \larrow \R^N$ be an immersion in $C^2(\Om)^N$ which solves $Du \ot Du :D^2u=0$. Then, $u$ is a Rank-One Absolute Minimal on $\Om$.
\end{lemma}

\BPL \ref{l1a}. If $\rk(Du)=n\leq N$ and $Du \ot Du :D^2u=0$ on $\Om$, then 
\beq \label{2.12b}
DuD\left(\frac{1}{2}|Du|^2\right)\ = \ 0.
\eeq
For each $x\in \Om$, the linear map $Du(x) : \R^n \larrow \R^N$ is injective and as such there exists a left inverse $(Du(x))^{-1}$. Hence, we obtain
\beq \label{2.12c}
(Du)^{-1}DuD\left(\frac{1}{2}|Du|^2\right)\ = \ 0
\eeq
which implies $D\big(\frac{1}{2}|Du|^2\big)=0$. Consequently, $u$ is a solution of the Eikonal equation since $|Du|$ is constant on connected components of $\Om$. Fix $D\subset \subset \Om$, $g\in C^1_0(D)$ and $\xi \in \mS^{N-1}$. We may assume $D$ is connected. Then, since $g|_{\p D}\equiv 0$, there exists an interior critical point $\bar{x}\in D$ of $g$. By using that $Dg(\bar{x})=0$, we estimate
\begin{align}
\big\|D(u+g\xi)\big\|_{L^\infty(D)}\ &= \ \sup_D \big|Du \, +\, \xi \ot Dg\big| \nonumber\\
            & \geq \ \big|Du(\bar{x})\, +\, \xi\ot Dg(\bar{x})\big|\\
            &=\ |Du(\bar{x})|  \nonumber\\
&=\ \big\|Du\big\|_{L^\infty(D)}. \nonumber
\end{align}
The lemma follows.
\qed

We have not been able to verify the validity of Lemma \ref{l1a} in the case of submersions for $N<n$, but we believe it is true. The difficulty lies in that the functionals $\Lip(u,\Om)$ and $\ess\,\sup_\Om |Du|$ are equivalent but not equal and standard scalar arguments as in \cite{A3} fail (cf.\ \cite{SS}).

\bl \label{l2} Let $u : \Om \sub \R^n \larrow \R^N$ be an immersion in $C^2(\Om)^N$ with $\infty$-Minimal area.Then, $u$ is normally $\infty$-Harmonic and solves $|Du|^2[Du]^\bot\De u=0$ on $\Om$.
\el

\BPL \ref{l2}. Fix $x\in \Om$, $0<\e<\dist(x,\p \Om)$ and $0<\de<1$. Fix also a normal vector field $\nu \in \Gamma\big([Du]^\bot \big)$ and an $h \in C^1\big(\overline{\mB_\e(x)} \big)$. We may assume that $\nu$ is a unit vector field. By differentiating the equation $|\nu|^2=1$ we obtain
\beq \label{2.12}
\nu^\top D\nu \ =\ 0.
\eeq
Moreover, by differentiating $\nu^\top Du = 0$ we obtain
\beq \label{2.13}
D\nu^\top Du\ =  \ -\nu^\top D^2u
\eeq
and by contracting \eqref{2.13} we get
\beq \label{2.14}
D\nu:Du\ =  \ -\nu^\top \De u.
\eeq
We set $w:=u+\de h\nu$. Then, we use that $\nu^\top$ annihilates $Du,D\nu$ and calculate:
\begin{align} \label{2.15}
|Dw|^2\ &=\ \big| Du \ + \ \de(\nu\ot Dh\ + \ hD\nu)\big|^2 \nonumber\\
             &=\ \big| (Du \ + \ \de hD\nu) \ + \ \de \nu\ot Dh\big|^2\\
             &=\ \big| Du \ + \ \de hD\nu \big|^2\ + \ \de^2|\nu|^2 |Dh|^2  \nonumber\\
             &=\ | Du|^2 \ + \ 2\de h\big(D\nu:Du \big)\ + \ \de^2\big(h^2|D\nu|^2+ |Dh|^2\big). \nonumber
\end{align}
By  \eqref{2.3} and equations \eqref{2.14}, \eqref{2.15}, we have
\begin{align} \label{2.16}
\|Du\|_{L^\infty(\mB_\e(x))}^2\ &\leq \ \|Dw\|_{L^\infty(\mB_\e(x))}^2 \nonumber\\
             &\leq \|Du\|_{L^\infty(\mB_\e(x))}^2\ +\ 2\de\sup_{\mB_\e(x)}\big\{h(D\nu:Du)\big\} +\ O(\de^2)\\
             & = \ \|Du\|_{L^\infty(\mB_\e(x))}^2\ -\ 2\de\min_{\overline{\mB_\e(x)}}\big\{h (\nu^\top \De u)\big\}\ +\ O(\de^2). \nonumber
           \end{align}
Hence, as $\de \ri 0$ we obtain
\beq \label{2.17}
\min_{\overline{\mB_\e(x)}}\big\{h (\nu^\top \De u)\big\}\ \leq \ 0.
\eeq
We now choose as $h$ the constant function $h:= \sgn \big( (\nu^\top \De u)(x)\big)$ and by \eqref{2.17} as $\e \ri 0$ we get $|(\nu^\top \De u)(x)|=0$. Since $\nu$ is an arbitrary unit normal vector field and $x$ is an arbitrary point, we get $[Du]^\bot\De u=0$ on $\Om$ and the lemma follows.                              \qed

\begin{remark} Equation \eqref{2.13} expresses the shape operator in the normal direction $\nu$ in terms of the second fundamental form of the submanifold $u(\Om)$. 
\end{remark}

\bl \label{l3} Let $u :\Om \sub \R^n \larrow \R^N$ be an immersion in $C^2(\Om)^N$. Suppose $u$ solves $|Du|^2[Du]^\bot\De u=0$ on $\Om$. Then, \emph{for all} $p \in [2,\infty]$, $u(\Om)$ has $p$-Minimal Area:
\beq \label{2.18}
\left.
\begin{array}{l}
D \subset \subset \Om, \\
h\in C^1(\bar{D}), \\
\nu \in \Gamma([Du]^\bot)
\end{array}
\right\} \ \ \Longrightarrow \ \
\big\|Du\big\|_{L^p(D)}\ \leq \ \big\|D(u+h\nu)\big\|_{L^p(D)}.
\eeq
Conversely, if \emph{for some}  $p\in [2,\infty)$ the image $u(\Om)$ has $p$-Minimal area, then $u$ solves $|Du|^2[Du]^\bot\De u=0$ on $\Om$.
\el

\BPL \ref{l3}. We begin with two differential identities. For any unit vector field $\nu \in \Gamma \big([Du]^\bot\big)$, $D\subset \subset \Om$, $h\in C^1(\bar{D})$, $\e \in \R$ and $p\geq 2$ we have 
\begin{align} \label{2.21a}
\frac{d}{d \e}\int_D \big|D(u+\e h&\nu)\big|^{p}\ = \ p\int_D\big|D(u+\e h\nu)\big|^{p-2}D(u+\e h\nu):D(h\nu),\\
\frac{d^2}{d \e^2} \int_D \big|D(u+\e h&\nu)\big|^{p}\ = \ p\int_D 
\big|D(u+\e h\nu)\big|^{p-2} \big| D(h\nu)\big|^2 \label{2.22a}\\
&+\ p(p-2)\int_D \big|D(u+\e h\nu)\big|^{p-4}\Big(D(u+\e h\nu): D(h\nu)\Big)^2.  \nonumber
\end{align}
Evidently, the function $\e \mapsto \int_D \big|D(u+\e h\nu)\big|^{p}-\int_D|Du|^p$  vanishes at $\e=0$ and by \eqref{2.22a} it is convex. By \eqref{2.12}, \eqref{2.14} and \eqref{2.21a} we have
\begin{align} \label{2.23a}
\frac{d}{d \e}\Big|_{\e=0} \int_D \big|D(u+\e h\nu)\big|^{p}\ &=\ p\int_D |Du|^{p-2}Du:D(h\nu) \nonumber\\
&= \ p\int_D |Du|^{p-2}Du:\big( hD\nu\ +\ \nu\ot Dh\big)\\
&= \ p\int_D |Du|^{p-2}(Du : D\nu)h \nonumber\\
&= \ -p\int_D |Du|^{p-2}(\nu^\top \De u)h.  \nonumber
\end{align} 
Since $|Du|>0$ on $\Om$ and $\nu$, $h$ are arbitrary, by \eqref{2.23a} we have  $|Du|^2[Du]^\bot \De u=0$ on $\Om$ if and only if
\beq \label{2.24a}
\int_D \big|Du\big|^{p} \ \leq \int_D \big|D(u+\e h\nu)\big|^{p}
\eeq
which means that $u(\Om)$ has $p$-Minimal area. By rescaling \eqref{2.24a} and letting $p\ri \infty$ we see that $u(\Om)$ has $\infty$-Minimal area as well when $u$ solves $|Du|^2[Du]^\bot \De u=0$ on $\Om$. The lemma has been established.            \qed

\ms
In view of Lemmas \ref{l1}, \ref{l1a}, \ref{l2} and \ref{l3}, Theorem \ref{1} follows. 

\begin{remark} Actually, in Lemma \ref{l3} we proved the stronger statement that the normal system $[Du]^\bot \De u=0$ characterizes immersions whose image $u(\Om)$ has $p$-Minimal area for \emph{any $p\in [2,\infty]$} and not only $p=\infty$.
\end{remark}

\subsection{Maximum and Minimum Principles for $|Du|$ for $\infty$-Harmonic Maps} We conclude this section by establising maximum and minimum principles for the gradient of $\infty$-Harmonic maps of full rank.

\begin{proposition}[Gradient Maximum-Minimum Principles] Suppose \label{pr1} $u :\Om \sub \R^n \larrow \R^N$ is in $C^2(\Om)^N$, $\infty$-Harmonic and of full rank. Then, for any $D \subset \subset \Om$ we have:
 \begin{align} 
\sup_{D}|Du|\ \leq \ \max_{\p D}|Du|, \label{2.23}\\
\inf_{D}|Du|\ \geq \ \min_{\p D}|Du|. \label{2.24}
\end{align}
\end{proposition}
In the case of submersions, the proof follows closely the ideas of Aronsson in [\cite{A3}, p.\ 558] and relates to the arguments of Capogna and Raich in [\cite{CR}, th.\ 1.1] perfomed for the special case of diffeomorphisms but for a different Hamiltonian in place of the Euclidean norm. The proof is based on the usage of the following improved modification of the gradient flow with parameters introduced in \cite{K2}:
\begin{lemma} \label{l5}
Let $u :\Om \sub \R^n \larrow \R^N$ be in $C^2(\Om)^N$. Consider the gradient flow
\beq \label{2.26}
\left\{
\begin{array}{l}
\dot{\ga}(t)\ = \ \left(\dfrac{|Du|^2}{|\xi^\top Du|^2}\xi^\top Du\right)\big(\ga(t)\big), \ \ t\neq 0,\ms\\
\ga(0)\ = \ x,
\end{array}
\right.
\eeq
for $x\in \Om$, $\xi \in \mS^{N-1}\set [Du(x)]^\bot$. Then, we have the differential identities
\begin{align}
\ \ \ \ \ \frac{d}{d t} \Big(\dfrac{1}{2}\big|Du\big(\ga(t)\big)\big|^2\Big)\ &= \ \left(\frac{|Du|^2}{|\xi^\top Du|^2}\xi^\top Du \ot Du :D^2 u \right)\big(\ga(t)\big), \label{2.27}\\
\frac{d}{d t}\Big(\xi^\top u\big(\ga(t)\big)\Big)\ &= \ \big| Du \big(\ga(t)\big)\big|^2, \label{2.28}
\end{align}
which imply $Du \ot Du:D^2u=0$ on $\Om$ if and only if $|Du\big(\ga(t)\big)|$ is constant along trajectories $\ga$ and $t\mapsto \xi^\top u\big(\ga(t)\big)$ is \emph{affine}.
\end{lemma}

We refrain from presenting the elementary proof of Lemma \ref{l5} which follows by simple calculations. We observe that in the scalar case of $N=1$, we have $\xi \in \{-1,+1\}$ and \eqref{2.26} reduces to the well known gradient flow (\cite{C}).

\BPP \ref{pr1}. 
Consider first the case of immersions where $\rk(Du)=n\leq N$. By arguing as in Lemma \ref{l1a}, it follows that $|Du|$ is constant on connected components of $\Om$. Hence, \eqref{2.23} and \eqref{2.24} follow.

For the case of submersions where $\rk(Du)=N\leq n$, fix $D \subset \subset \Om$, $x\in D$ and $\xi \in \mS^{N-1}$ and consider the gradient flow \eqref{2.26}. Since $\rk(Du)=N\leq n$, for each $y\in \Om$ the linear map $Du(y)^\top : \R^N \larrow \R^n$ is injective and hence $|\xi^\top Du|>0$ on $\Om$. Hence, the flow is globally defined on $\Om$ for all parameters $\xi$. By \eqref{2.27}, $|Du\big(\ga(t)\big)|=|Du(x)|$ and by \eqref{2.28} the trajectory $\ga$ reaches $\p D$ in finite time since $D$ is bounded while
\begin{align}
\xi^\top u( \ga(t))\ - \ \xi^\top u(x)\ &=\, t|Du(x)|^2.
\end{align}
Hence, there exists $t^+(x)>0$ such that $\ga(t^+(x)) \in \p D$. Consequently,
\begin{align}
\sup_{D}|Du| \ &= \ \sup_{x\in D}|Du(x)| \nonumber\\
   &= \ \sup_{x\in D}|Du\big(\ga(t^+(x))\big)|\\
& \leq \  \max_{\p D} |Du|\nonumber
\end{align}
and similarly we obtain $\inf_D|Du|\geq \min_{\p D}|Du|$. The proposition follows.                      \qed

\section{Connections to Minimal Surfaces.} \label{section3}

In this section we restrict attention to $2$-dimensional $\infty$-Harmonic immersions $u:\Om \sub \R^2 \larrow \R^3$ and draw tight connections to Differential Geometry. We show that abstract smooth minimal surfaces of $\R^3$ can be characterized as those that can be locally parameterizable by isothermal immersions which are normally $\infty$-Harmonic, that is by conformal coordinate maps with $\infty$-minimal area (Definition \ref{def1}). Moreover, we show that isothermal $\infty$-Harmonic maps are rigid and they always have affine range.

We begin with two differential identites which connect $\De_\infty$ to the geometry of the range of conformal $\infty$-Harmonic maps. Interestingly, the lemma holds for conformal maps \emph{with degeneracies}, that is when there exists $f\in C^1(\Om)$ such that $Du^\top Du =f^2 I$ on $\Om$ but $f$ may have \emph{zeros}.

\bl \label{l4} Let $u : \Om \sub \R^2 \larrow \R^3$ be a conformal map in $C^2(\Om)^3$. Then, we have the identities
\begin{align}
\frac{1}{2}|Du|^2\ &=  \ \sqrt{\det(g)}, \label{3.1}\\
|Du|^2[Du]^\bot \De u \ &= \ 4\det(g)\H \label{3.2},
\end{align}
where $H$ is the mean curvature vector of $u\big(\Om \set \{\det(g)=0\}\big)$ and $g=Du^\top Du$ is the  induced Riemannian metric, that is $\sqrt{\det(g)}$ equals the Jacobian $Ju$.
\el

\begin{remark} \eqref{3.1} is valid also for maps $u : \Om \sub \R^2 \larrow \R^N$ for any $N\geq 3$. 
\end{remark}

\BPL \ref{l4}. By assumption there is $f\in C^1(\Om)$ such that $Du^\top Du =f^2 I$ on $\Om$. If $D_xu$, $D_y u$ $: \Om \sub \R^2 \larrow \R^3$ denote the two partial derivatives of $u$, then we have $|D_xu|^2=|D_y u|^2=f^2$ and $D_x u^\top D_y u=0$. Moreover, $f^2=\frac{1}{2}\tr(Du^\top Du)=\frac{1}{2}|Du|^2$. Hence, we have
\begin{align} \label{3.3}
|Du|^2\ &=\ |D_x u|^2\ + \ |D_y u|^2 \nonumber\\
&=\ 2\, f^2  \nonumber\\
&=\ 2\big(|D_x u|^2|D_y u|^2\big)^{\frac{1}{2}}\\
&=\ 2\big(|D_x u|^2|D_y u|^2\, -\, (D_x u^\top D_y u )^2\big)^{\frac{1}{2}} \nonumber\\
&=\ 2\det\big(Du^\top Du\big)^{\frac{1}{2}}. \nonumber
\end{align}
Hence, \eqref{3.1} follows. Let now $\nu$ be the normal vector field over $\Om \set \{Ju\neq 0\}$. Then, the mean curvature vector $\H$ of the immersion $u : \Om \set \{Ju\neq 0\} \larrow \R^3$ is
\beq \label{3.4}
\H\ = \ \frac{ |D_x u|^2({\nu}^\top D^2_{yy}u) \, + \, |D_y u|^2({\nu}^\top D^2_{xx}u)  \, -\,2 (D_xu^\top D_y u)({\nu}^\top D^2_{xy}u)  }{2\left(|D_x u|^2|D_y u|^2\, -\, (D_x u^\top D_y u )^2\right)}\nu.
\eeq
Since $u$ is conformal, we have
\begin{align}
\H\ &= \ \frac{f^2 \big({\nu}^\top D^2_{yy}u \, + \, {\nu}^\top D^2_{xx}u  \big)}{2\left(|D_x u|^2|D_y u|^2\, -\, (D_x u^\top D_y u )^2\right)}\nu \nonumber
\\
& = \ \frac{f^2}{2\det\big(Du^\top Du\big)}{(\nu\ot \nu) \big(D^2_{yy}u \, + \, D^2_{xx}u\big) }\\
&=\ \frac{1}{2\det\big(Du^\top Du\big)^{\frac{1}{2}} } [Du]^\bot \De u. \nonumber
\end{align}
Hence, on $\Om \set \{\det(g)\neq 0\} $ we have
\beq \label{3.6}
2\sqrt{\det(g)}\H\ = \ [Du]^\bot \De u.
\eeq
Equation \eqref{3.6} readily leads to \eqref{3.2} on $\Om \set \{\det(g)\neq 0\} $ and extends to $\Om$ since both sides vanish on $\{Ju=0\}$. The lemma follows.                                 \qed

\ms

Formulas \eqref{3.1} and \eqref{3.4} readily lead to the next

\bcor \label{cor2} Let $u : \Om \sub \R^2 \larrow \R^3$ be an immersion in $C^2(\Om)^3$.

\ms (i) If $u$ is conformal, the surface area of $u(\Om)$ is
\beq
\mathcal{H}^2\big(u(\Om)\big)\ = \ \int_\Om \frac{1}{2}|Du|^2,
\eeq
where $\mathcal{H}^2$ is the $2$-dimensional Hausdorff measure.

\ms (ii) If  $|D_x u|=|D_y u|$ and $|Du|^2[Du]^\bot\De u=0$  on $\Om$, the mean curvature vector of $u(\Om)$ is given by
\beq
\H\ = \ -\frac{1}{(Ju)^2}\big(D_x u ^\top D_y u\big)[Du]^\bot D^2_{xy} u.
\eeq
Hence, $u(\Om)$ is minimal if and only if either $D^2_{xy}u$ is tangential or $u$ is conformal.
\ecor

\begin{theorem}[Minimal Surfaces and conformal $\infty$-Harmonic maps] \label{th2} Let $\S \sub \R^3$ be a $C^2$ surface, with the induced Riemannian metric. Then, 
\ms

(i) $\S$ is minimal if and only if $\S$ has an atlas of isothermal normally $\infty$-Harmonic parametrizations.
\ms

(ii) $\S$ is locally isometric to $(\R^2,c^2I)$ for some $c\in \R$ if and only if $\S$ has an atlas of isothermal tangentially $\infty$-Harmonic parametrizations.

\ms

(iii) $\S$ is contained in an affine plane of $\R^3$ if and only if $\S$ has an atlas of isothermal $\infty$-Harmonic parametrizations.
\end{theorem}

Theorem \ref{th2} readily implies the following

\bcor[Rigidity of conformal $2$-dimensional $\infty$-Harmonic maps] \label{cor3} If $u : \Om \sub \R^2 \larrow \R^3$ is in $C^2(\Om)^3$, conformal and $\infty$-Harmonic, then $u(\Om)$ is contained into an affine plane of $\R^3$.
\ecor

\BPT \ref{th2}. We begin by recalling the standard fact that every point $p\in \S$ of a smooth surface has an isothermal parametrization $u : \Om \sub \R^2 \larrow \S$, where $p=u(x_0,y_0)$ for some $(x_0,y_0)^\top \in \Om$. For every such $u$, there exists an $f\in C^1(\Om)$ such that $Du^\top Du=f^2I$ on $\Om$.
Then, (i) follows by observing that identity \eqref{3.2} implies that $[Du]^\bot\De u=0$ if and only if the mean curvature of $\S$ vanishes. 

To see (ii), first assume that $Du\ot Du:D^2u=0$ on $\Om$. Since $u$ is an immersion, equations \eqref{2.12b} and \eqref{2.12c} imply that $\frac{1}{2}|Du|^2=c^2$ for some $c\in \R$. By the following elementary identities which are valid for $2$-dimensional conformal maps
\begin{align}
\tr\big(Du^\top Du-c^2I\big)\ &= \ |Du|^2 \ -\ 2c^2 \label{3.8}\\
\det\big(Du^\top Du-c^2I\big)\ &= \ \Big(\frac{1}{2}|Du|^2 \ -\ c^2 \Big)^2
\end{align}
we obtain that $Du^\top Du=c^2I$. Hence, $\S$ is locally isometric to $(\R^2,c^2I)$. Conversely, if $\S$  is locally isometric to $(\R^2,c^2I)$, then we have that $Du^\top Du=c^2I$ and hence by \eqref{3.8} we have $|Du|^2=2c^2$, which implies  $Du\ot Du:D^2u=0$ on $\Om$.

Finally, $(iii)$ follows by observing that the only minimal surfaces which are locally isometric to $(\R^2,c^2I)$ are portions of affine planes of $\R^3$. Indeed, fix an isothermal parametrization  $u : \Om \sub \R^2 \larrow \S$ of the surface $\S$. Then, if $\S$ is minimal and isometric to $(\R^2,c^2I)$, both the principal curvatures vanish since the mean curvature and the Gauss curvature vanish. Hence, the shape operator vanishes and as such $u(\Om)$ is contained into an affine plane. The converse implication is obvious.                                \qed

\begin{remark} The results of this paper extend with little extra cost to general supremal functionals $\ess\, \sup_\Om H(Du)$ for a convex Hamiltonian $H \in C^2(\R^N \ot \R^n)$ and the respective Aronsson system studied in \cite{K2,K3}. We just observe that \eqref{2.14} generalizes to
\beq \label{3.11}
D \nu : H_P(Du)\ = \ -\nu^\top H_{PP}(Du):D^2u
\eeq
and \eqref{3.11} follows by differentiating the equation $\nu^\top H_P(Du)=0$. The latter says that $\nu$ is a section of the vector bundle over $u(\Om)$ with fibers $[H_P(Du(x))]^\bot$ where $x\in \Om \sub \R^n$.
\end{remark}

In the forthcoming work \cite{K5} we present a theory of non-differentiable solutions which applies to fully nonlinear PDE systems and extends Viscosity Solutions to the general vector case. This approach is based on the existence of an extremality principle which applies to maps. In this context, we consider  the existence of solution to the Dirichlet problem for \eqref{1.1}.

\ms

\noi \textbf{Acknowledgement.} The observation in the proof of $(iii)$ of Theorem \ref{th2} is due to S.\ Aretakis. I thank Y.\ Yu, J.\ Manfredi, L.C.\ Evans and L.\ Capogna for their interest and encouragenment. I am indebted to the anonymous referee for his valuable comments which improved both the content and the appearance of the paper.

\bibliographystyle{amsplain}

\end{document}